\documentclass[12pt]{amsart}
\usepackage{latexsym}
\usepackage{amssymb,times,fullpage}
\theoremstyle{plain}
\newtheorem{lemma}{Lemma}[section]
\newtheorem{proposition}[lemma]{Proposition}
\newtheorem{theorem}[lemma]{Theorem}

\theoremstyle{definition}
\newtheorem{definition}[lemma]{Definition}
\newtheorem{remark}[lemma]{Remark}
\newtheorem{example}[lemma]{Example}

\newcommand{\R}{\mathbb{ R}}

\newcommand{\Z}{\mathbb{ Z}}

\def\M1g{{\mathcal M}_{g,1}}
\def\I1g{{\mathcal I}_{g,1}}

\newcommand\Aut{\operatorname{Aut}}

\newcommand\Out{\operatorname{Out}}

\newcommand\Cal{\operatorname{Cal}}
\newcommand\Ker{\operatorname{Ker}}

\newcommand\Diff{\operatorname{Diff}}

\title[Stable length in stable groups]{Stable length in stable groups\footnote{\copyright \ D.~Kotschick 2006--2008}}
\author[D.~Kotschick]{D.~Kotschick}
\address{Mathematisches Institut, Ludwig-Maximilians-Universit\"at M\"unchen,
Theresienstr.~39, 80333 M\"unchen, Germany}
\email{dieter@member.ams.org}
\date{June 16, 2008; MSC 2000: primary 20F69, secondary 20F12, 57M07}

\begin{document}

\begin{abstract}
We show that the stable commutator length vanishes for certain groups defined as infinite unions of smaller 
groups. The argument uses a group-theoretic analogue of the Mazur swindle, and goes back to the works
of Anderson, Fisher, and Mather on homeomorphism groups.
\end{abstract}


\maketitle

\section{Introduction}

In this paper we show that the stable commutator length vanishes for certain groups defined as unions of 
subgroups that have many conjugate embeddings that commute element-wise. The argument used is a 
group-theoretic analogue of the Mazur swindle. Informally, it can be paraphrased by saying that as the size of a table grows, 
it becomes easier to sort objects on its surface, so that commutation can be done very efficiently on a table of 
infinite size\footnote{I owe this description to N.~A'Campo.}. In the words of Po\'enaru~\cite{Po}: ``The infinite 
comes with magic and power.''

Let $\Gamma$ be a group. The commutator length $c(g)$ of an element $g$ in the commutator subgroup 
$[\Gamma,\Gamma]\subset\Gamma$ is the minimal number of factors needed to write $g$ as a product of elements that 
can be expressed as single commutators. The stable commutator length of $g$ is defined to be 
$$
\vert\vert g\vert\vert =\lim_{n\rightarrow\infty}\frac{c(g^n)}{n} \ .
$$
Direct proofs of the vanishing of the stable commutator length sometimes proceed by showing a lot more, namely
that the commutator length itself is bounded. A prototypical argument for this in the context of homeomorphism groups 
goes back to Anderson~\cite{A} and Fisher~\cite{F}, and was later used and refined by Mather~\cite{M1} and by
Matsumoto--Morita~\cite{MM}. This argument uses an infinite iteration that leads to complicated behaviour near
one point, and is therefore not suitable for the study of diffeomorphism groups. The argument we employ here is 
a variation on this classical one. It works for diffeomorphism groups because instead of an infinite iteration we do 
only a finite iteration, however it is important that the finite number of iterations can be taken to be arbitrarily large. 
This argument does not prove that the 
commutator length is bounded, it only proves the weaker conclusion that the stable commutator length vanishes.

Bavard~\cite{B}, using a Hahn--Banach argument influenced by that of Matsumoto and Morita~\cite{MM}, 
showed that the vanishing of the stable commutator length on the commutator 
subgroup is equivalent to the injectivity of the comparison map $H^2_b(\Gamma;\R)\rightarrow H^2(\Gamma;\R)$,
where $H^*_b(\Gamma)$ denotes bounded group cohomology in the sense of Gromov. Another way to express
this condition is to say that every homogeneous quasi-homomorphism $\Gamma\rightarrow\R$ is in fact a 
homomorphism. If $\Gamma$ is a perfect group, then the stable commutator length is defined on the whole
of $\Gamma$, and vanishes identically if and only if there are no non-trivial homogeneous quasi-homomorphisms
from $\Gamma$ to $\R$.

Bavard's result suggests that the commutator calculus arguments discussed above should have a dual formulation
in terms of quasi-ho\-mo\-mor\-phisms. It is this dual formulation that we shall discuss, although we could equally well
argue directly on the commutator side. Thus, what we prove for certain groups is that every homogeneous quasi-homomorphism
is a homomorphism. Using Bavard's result this is equivalent to the vanishing of the stable commutator length.

In Section~\ref{s:alg} we give the algebraic mechanism behind the vanishing results we shall prove.
Section~\ref{s:discrete} applies this mechanism to the stable mapping class group, the braid group on
infinitely many strands, and the stable automorphism groups of free groups. In Section~\ref{s:diff} we 
give applications to diffeomorphism groups, and in Section~\ref{s:BIP} we compare our methods and results
to those of Burago, Ivanov and Polterovich~\cite{BIP}.

We refer the reader to~\cite{B,qmorph} for background on quasi-ho\-mo\-mor\-phisms. In fact, the proofs of the vanishing 
results below are somewhat reminiscent of the discussion of weak bounded generation in~\cite{qmorph}.

\section{The algebraic vanishing result}\label{s:alg}

A map $f\colon\Gamma\rightarrow\R$ is called a quasi-homomorphism 
if its deviation from being a homomorphism is bounded; in other words, 
there exists a constant $D(f)$, called the defect of $f$, such that
$$
\vert f(xy)-f(x)-f(y)\vert\leq D(f)
$$
for all $x,y\in\Gamma$. We will always take $D(f)$ to be the smallest number 
with this property, i.e.~it is the supremum of the left hand sides over all $x$
and $y\in\Gamma$.

Every quasi-homomorphism can be homogenized by defining
$$
\varphi (g) = \lim_{n\rightarrow\infty}\frac{f(g^{n})}{n} \ .
$$
Then $\varphi$ is again a quasi-homomorphism, is homogeneous in the 
sense that $\varphi(g^{n})=n\varphi(g)$, and is constant on conjugacy 
classes. (Compare~\cite{B}, Proposition~3.3.1.) Throughout this paper we shall only consider 
homogeneous quasi-homomorphisms.

We begin with the following preliminary result.
\begin{lemma}\label{l:pre}
Let $\varphi\colon\Gamma\rightarrow\R$ be a homogeneous quasi-homomorphism. Then the following two 
properties hold:
\begin{enumerate}
\item If $x$, $y\in\Gamma$ commute, then $\varphi (xy)=\varphi (x)+\varphi (y)$.
\item If $\varphi$ vanishes on every single commutator, then $\varphi$ is a homomorphism.
\end{enumerate}
\end{lemma}
\begin{proof}
By homogeneity we have the following:
\begin{alignat*}{1}
\vert \varphi(xy)-\varphi(x)-\varphi(y)\vert &= \lim_{n\rightarrow\infty}\frac{1}{n}\vert\varphi ((xy)^n)+\varphi (x^{-n})+
\varphi (y^{-n})\vert\\
&\leq \lim_{n\rightarrow\infty}\frac{1}{n}\left(\vert \varphi((xy)^nx^{-n}y^{-n})\vert+2D(\varphi)\right)\\
&=\lim_{n\rightarrow\infty}\frac{1}{n}\vert \varphi((xy)^nx^{-n}y^{-n})\vert \ .
\end{alignat*}
If $x$ and $y$ commute, then $\varphi((xy)^nx^{-n}y^{-n})=0$ for every $n$, giving the first statement.

Assume now that $\varphi$ vanishes on single commutators. As $(xy)^nx^{-n}y^{-n}$ can be expressed as the product of
$\frac{n}{2}+c$ commutators, see~\cite{B}, the right hand side of the formula is bounded above by $\frac{1}{2}D(\varphi )$.
However, taking the supremum of the left hand side over all $x$ and $y$, we get the defect $D(\varphi )$. Thus
$D(\varphi )\leq \frac{1}{2}D(\varphi )$, showing that the defect vanishes, and $\varphi$ is a homomorphism.
\end{proof}

Here is the main mechanism for the vanishing theorems.
\begin{proposition}\label{p:hom}
Let $\Lambda\subset\Gamma$ be a subgroup with the property that there is an arbitrarily large number of conjugate embeddings 
$\Lambda_i\subset\Gamma$ of $\Lambda$ in $\Gamma$ with the property that elements of $\Lambda_i$ and of $\Lambda_j$ commute
with each other in $\Gamma$ whenever $i\neq j$. Then every homogeneous quasi-homomorphism on $\Gamma$ restricts to $\Lambda$ 
as a homomorphism.
\end{proposition}
Note that homogeneous quasi-homomorphisms are constant on conjugacy classes, so that the restriction of $\varphi$ to $\Lambda_i$
is independent of $i$.
\begin{proof}
By the second part of Lemma~\ref{l:pre} we only have to prove that $\varphi([x,y])=0$ for any $x,y\in\Lambda$. Let $x_i$, $y_i\in\Gamma$
be the images of $x$ and $y$ under the embedding $\Lambda_i\subset\Gamma$. 
We have the following equalities:
\begin{alignat*}{1}
n\varphi([x,y]) &=\varphi([x_1,y_1])+\ldots +\varphi([x_n,y_n]) \\
&=\varphi([x_1,y_1]\ldots [x_n,y_n]) \\
&=\varphi([x_1\ldots x_n,y_1\ldots y_n]) \ ,
\end{alignat*}
where the first one comes from the constancy of $\varphi$ on conjugacy classes, the second is the first part of Lemma~\ref{l:pre} applied to the 
commuting embeddings, and the third follows directly from the commuting property of the embeddings $\Lambda_i\subset\Gamma$.
On the right hand side $\varphi$ is applied to a single commutator in $\Gamma$, and therefore the right hand side is bounded in absolute
value by the defect of $\varphi$ on $\Gamma$. However, if $\varphi([x,y])\neq 0$, then the left hand side is unbounded because by 
assumption we can make $n$ arbitrarily large. Thus $\varphi([x,y])=0$ for all $x,y\in\Lambda$.
\end{proof}

There are several ways in which this mechanism can be  applied to obtain the vanishing of the  stable commutator length in 
various groups. We try to give a general statement, in the hope of unifying several different applications of the argument.
\begin{theorem}\label{t:alg}
Let $\Gamma$ be a group in which every element can be decomposed as a product of some fixed number $k$ of elements contained
in distinguished subgroups $\Lambda\subset\Gamma$. If each $\Lambda$ is perfect and has the property in Proposition~\ref{p:hom}, then
the stable commutator length of $\Gamma$ vanishes.
\end{theorem}
\begin{proof}
As the subgroups $\Lambda$ are assumed perfect, the restrictions of quasi-homomorphisms on $\Gamma$ to $\Lambda$ vanish,
because by Proposition~\ref{p:hom} they are homomorphisms. Therefore the value of a quasi-homomorphism on every element of $\Gamma$
is bounded by $k-1$ times the defect. But every bounded homogeneous quasi-homomorphism is trivial.
\end{proof}

\section{Applications to discrete groups}\label{s:discrete}

\subsection{The stable mapping class group}

Let $\Gamma_g^1$ be the group of isotopy classes of diffeomorphisms with compact support in the interior of a compact surface $\Sigma_g^1$
of genus $g$ with one boundary component. Attaching a two-holed torus along the boundary defines the stabilization homomorphism 
$\Gamma_g^1\rightarrow\Gamma_{g+1}^1$. The stable mapping class group $\Gamma_{\infty}$ is defined as the limit
$$
\Gamma_{\infty} = \lim_{\stackrel{\longrightarrow}{g}}\Gamma_g^1 \ .
$$
For $g\geq 3$ the groups $\Gamma_g^1$ are perfect, see~\cite{Powell}, hence $\Gamma_{\infty}$ is also perfect.
Recall that in~\cite{EK} it was proved that the stable commutator length is non-trivial on every $\Gamma_g^1$ with $g\geq 2$,
see also~\cite{BF,BK,CF,qmorph}. In contrast with this we have:
\begin{theorem}\label{t:sMCG}
The stable commutator length for $\Gamma_{\infty}$ vanishes identically.
\end{theorem}
\begin{proof}
We apply Theorem~\ref{t:alg} with $k=1$. The subgroups $\Lambda$ are the images of the $\Gamma_g^1$. In detail, every 
element of $\Gamma_{\infty}$ is in the image of some $\Gamma_g^1=\Lambda$. This has arbitrarily large numbers of commuting 
conjugate embeddings in $\Gamma_{\infty}$ given by taking the boundary connected sum of an arbitrarily large number of copies of 
the surface of genus $g$ with one boundary component. Any homogeneous quasi-homomorphism on $\Gamma_{\infty}$
restricts to (the image of) $\Gamma_g^1$ as a homomorphism. However, for $g\geq 3$ this group is perfect, and so the homomorphism 
vanishes.
\end{proof}
\begin{remark}
Theorem~\ref{t:sMCG} shows that the second bounded cohomology of mapping class groups does not stabilize. 
This contrasts sharply with the Harer stability theorem~\cite{Harer} for the ordinary group cohomology. 
\end{remark}
\begin{remark}
Theorem~\ref{t:sMCG} fits in nicely with the form of the estimates for the stable commutator length obtained in~\cite{EK,BK,qmorph}. 
The lower bounds given there for the stable commutator length of specific elements in $\Gamma_g^1$ go to zero for $g\rightarrow\infty$.
A similar phenomenon seems to appear in the work of Calegari and Fujiwara~\cite{CF}.
\end{remark}
\begin{remark}
The discussion in this subsection also applies to the stable mapping class groups for surfaces with several boundary components.
\end{remark}

\subsection{The braid group on infinitely many strands}

Let $B_n$ be the Artin braid group on $n$ strands. Adding strands defines injective stabilization homomorphisms
$B_n\longrightarrow B_{n+1}$. The braid  group on infinitely many strands is 
$$
B_{\infty} = \bigcup_n B_n \ .
$$
\begin{theorem}\label{t:braids}
Any homogeneous quasi-homomorphism on the infinite braid group $B_{\infty}$ is a homomorphism.
\end{theorem}
\begin{proof}
We apply Proposition~\ref{p:hom}. The subgroups $\Lambda$ are the $B_n$. Any two $x$, $y\in B_{\infty}$ 
are contained in some finite $B_n=\Lambda$. This has arbitrarily large numbers of commuting conjugate embeddings in 
$B_{\infty}$ given by considering braids on $n$ strands placed side by side. Any homogeneous quasi-homomorphism 
$\varphi$ on $B_{\infty}$ therefore restricts to $B_n$ as a homomorphism. Thus $\varphi (xy)=\varphi (x) +\varphi (y)$, 
showing that $\varphi$ is a homomorphism on $B_{\infty}$.
\end{proof}
We have chosen this formulation of the result because the braid groups have infinite Abelianizations, so that it does 
not make sense to speak of the stable commutator lengths of elements. The stabilization is compatible with Abelianization, 
so that $B_{\infty}$, like $B_n$ for finite $n$, has infinite cyclic Abelianization. The conclusion of Theorem~\ref{t:braids}
is that every homogeneous quasi-homomorphism is a constant multiple of the Abelianization homomorphism.

For finite $n$, the  braid groups do admit non-trivial homogeneous quasi-homomorphisms that are not proportional to 
the Abelianization map $B_n\longrightarrow\Z$, see~\cite{Baader,GG}. Thus Theorem~\ref{t:braids} shows that the bounded 
cohomology of braid groups does not stabilize. The quotients $B_n/C$ of finite braid groups modulo their centers have 
finite Abelianizations, so that all elements have powers that are products of commutators. Baader~\cite{Baader} proves 
some lower bounds for the stable commutator length of certain elements in $B_n/C$. These lower bounds tend to zero as $n\rightarrow\infty$.
This of course fits with Theorem~\ref{t:braids}. The corresponding homogeneous quasi-homomorphisms defined on $B_n/C$,
and, by composition, on $B_n$, do not extend to the braid group on infinitely many strands.

\subsection{Automorphism groups of free groups}

The canonical homomorphisms $F_n\longrightarrow F_n\star\Z=F_{n+1}$ give rise to injective homomorphisms
 $\Aut(F_n)\longrightarrow \Aut(F_{n+1})$. Thus we define
$$
\Aut_{\infty}(F) = \bigcup_{n}\Aut(F_n) \ .
$$
Note that this is smaller than $\Aut(F_{\infty})$.
\begin{theorem}\label{t:sAut}
The stable commutator length for $\Aut_{\infty}(F)$ vanishes identically.
\end{theorem}
\begin{proof}
The Abelianization of $\Aut(F_n)$ is of order $2$, and is stable. Thus, up to taking squares, all elements in $\Aut(F_n)$ 
and in $\Aut_{\infty}(F)$ are products of commutators.
We again apply Proposition~\ref{p:hom}. The subgroups $\Lambda$ are the $\Aut(F_n)$. In detail, every 
element of $\Aut_{\infty}(F)$ is contained in some $\Aut(F_n)=\Lambda$. This has arbitrarily large numbers of commuting 
conjugate embeddings in $\Aut_{\infty}(F)$ given by embedding $F_n\star F_n\star\ldots\star F_n$ in some large $F_N$. 
Although the elements in the different copies of $F_n$ do not commute, the induced embeddings of $\Aut(F_n)$ in
$\Aut_{\infty}(F)$ do commute and are conjugate to each other. Any homogeneous quasi-homomorphism on $\Aut_{\infty}(F)$
restricts to $\Aut(F_n)$ as a homomorphism. As the Abelianization of $\Aut(F_n)$ is finite, the homomorphism is trivial.
\end{proof}

One could also consider the outer automorphism groups $\Out(F_n)$, and of course quasi-ho\-mo\-mor\-phisms
on these induce quasi-homomorphisms on $\Aut(F_n)$ by composition with the projection. However, there is no natural
way  of stabilizing the outer automorphism groups.

Unlike for mapping class groups, no unbounded quasi-ho\-mo\-mor\-phisms are known on the automorphism groups of 
free groups. The groups $\Aut(F_n)$ are analogous to the so-called extended mapping class groups, consisting of the 
isotopy classes of all diffeomorphisms of surfaces. The usual mapping class groups we considered in Theorem~\ref{t:sMCG} 
are index $2$ subgroups of the extended mapping class groups. A single Dehn twist has non-zero stable commutator length 
in a mapping class group~\cite{EK,BK,qmorph}, but is conjugate to its inverse in the extended mapping class group~\cite{NAC}, 
so that in this latter group its stable commutator length must vanish. By analogy with this phenomenon, it may be more promising 
to look for quasi-homomorphisms on the special automorphism group $S\Aut(F_n)$, rather than
on $\Aut (F_n)$, where by  $S\Aut(F_n)$ we denote the automorphisms preserving an orientation on the Abelianization $\Z^n$
of $F_n$. Theorem~\ref{t:sAut} also applies to the  stabilized special automorphism group $S\Aut_{\infty}(F)$.

\section{Applications to diffeomorphism groups}\label{s:diff}

In this section we apply the algebraic vanishing result to some groups of diffeomorphisms. If $M$ is a smooth closed
manifold, we consider the identity component $G=\Diff_0(M)$ of the full diffeomorphism group. If $M$ is open, we consider
the identity component $G=\Diff_0^c(M)$ of the group of compactly supported diffeomorphisms. (Sometimes this notation
is redundant because the group of compactly supported diffeomorphisms may be connected.) In both cases we assume 
that the diffeomorphisms are of class $C^r$ with $1\leq r\leq\infty$ and $r\neq 1+\dim (M)$. The classical results of Herman, 
Thurston, Epstein and Mather then ensure that $G$ is a perfect group; see~\cite{Ban} and the references quoted there.

\subsection{The full diffeomorphism groups}

We would like to prove that the stable commutator length vanishes on the diffeomorphism groups. Unfortunately, we can only
achieve this goal in the following special cases:
\begin{theorem}\label{t:easy}
The stable commutator length vanishes for the following diffeomorphism groups:
\begin{enumerate}
\item $\Diff_0^c(D)$ for a ball $D$,
\item $\Diff_0^c(D\times M)$ for any $M$ and a ball $D$ of positive dimension,
\item $\Diff_0(S^n)$ for any sphere,
\item $\Diff_0(\Sigma^n)$ for any exotic sphere of dimension $\neq 4$.
\end{enumerate}
\end{theorem}
\begin{proof}
Let $D$ be an $n$-dimensional ball, and $G=\Diff_0^c(D)$ the group of diffeomorphisms of $D$ with compact support in the interior of $D$.
Fix an exhaustion of $D$ by smaller nested balls $D_i$ so that each $D_i$ has closure contained in the interior of $D_{i+1}$.
If $G_i$ is the group of diffeomorphisms of $D_i$ with compact support in the interior of $D_i$, then we have injective homomorphisms 
$G_i\rightarrow G_{i+1}$ induced by the  inclusion $D_i\subset D_{i+1}$, and 
$$
G = \bigcup_{i} G_i \ .
$$
We now apply Theorem~\ref{t:alg} with $k=1$. The subgroups $\Lambda$ are the $G_i$. In detail, every element of 
$G$ is contained in some $G_i$. This has arbitrarily large numbers of commuting conjugate embeddings in $G$, 
because we can embed arbitrarily large numbers of disjoint (smaller) balls in the annulus $D\setminus D_i$. Any 
homogeneous quasi-homomorphism on $G$ then restricts to $G_i$ as a homomorphism. But this group is perfect,
and so the homomorphism vanishes.
This completes the proof of the first claim.
The same argument works for $D\times M$ using the exhaustion by $D_i\times M$.

Now let $G=\Diff_0(S^n)$.
We can compose any $g\in G$ with a suitable $f_1^{-1}\in G$ supported in an open ball to achieve that $f_1^{-1} g$ has a fixed point. 
It follows that we can write $f_1^{-1} g=f_2 f_3$, with $f_3$ having support in a ball around the fixed point, and $f_2$ having support 
in the complement of a (smaller) ball around the fixed point. Because the manifold is a sphere, this complement is again a ball. 
Using $g=f_1 f_2 f_3$, we apply Theorem~\ref{t:alg} with $k=3$. Let $\Lambda=\Diff_c(D)$ be the group of compactly supported 
diffeomorphisms of a ball. This is perfect and admits infinitely many conjugate commuting embeddings in $G$.

The argument given for spheres also works for exotic spheres in dimensions $n\neq 4$, because they are all twisted spheres 
obtained from two standard balls by gluing along the boundary~\cite{Kos}.
\end{proof}

After this paper had been submitted, Tsuboi~\cite{T} proved a much stronger result. He proved that $\Diff_0(M)$ is uniformly
perfect for every closed manifold, with the possible exception of even-dimensional manifolds having no handle decomposition
without handles of middle dimension. The simplest $M$ for which the question is unresolved, is $T^2$. Because of Tsuboi's 
results, I have removed some material from the original version of this paper, which was concerned with the continuity, with
respect to the $C^0$ topology, of potential homogeneous quasi-homomorphisms on diffeomorphism groups. That discussion is 
rather delicate technically, and I hope to return to it elsewhere, unless it is rendered completely empty by a generalization of 
Tsuboi's work. 

\subsection{Diffeomorphism groups preserving a symplectic or volume form}

The first statement in Theorem~\ref{t:easy} can be phrased for the compactly supported diffeomorphism group of Euclidean
space $\R^n$ instead of a ball $D$. However, the two cases are rather different if we consider diffeomorphism groups 
preserving a symplectic or volume form $\omega$, and we fix our conventions so that $\R^n$ has infinite (symplectic) volume 
and a ball has finite volume.

Let $G=\Diff_0^c(\R^n,\omega)$ be the identity component of the group of compactly supported diffeomorphisms
preserving $\omega$. This is not perfect because the Calabi invariant 
\begin{equation}\label{eq:Cal}
\Cal\colon \Diff_0^c(\R^n,\omega)\longrightarrow\R
\end{equation}
is a non-trivial homomorphism. (See for example~\cite{Ban} for the definition of $\Cal$.)
By results of Thurston and Banyaga, see~\cite{Ban}, its kernel is a perfect group. 
It follows that $\Cal$ is the Abelianization homomorphism. So far there is no difference between the finite and 
infinite volume cases, and we have the same statements if we replace $(\R^n,\omega)$ by a ball of finite volume.
\begin{theorem}\label{t:Ham}
The stable commutator length vanishes on the kernel of the Calabi homomorphism in $\Diff_0^c(\R^n,\omega)$.
On $\Diff_0^c(\R^n,\omega)$ every homogeneous quasi-homomorphism is a constant multiple of the Calabi invariant.
\end{theorem}
\begin{proof}
As in the previous proof we can write the kernel of $\Cal$ as 
$$
G = \bigcup_{i} G_i \ ,
$$
with each $G_i$ consisting of those elements of $G$ supported in a ball of radius $i$, say. Then each element of $G$
is contained in some $G_i$, but this $G_i$ has arbitrarily large numbers of commuting conjugate embeddings given by
disjoint Hamiltonian displacements of $D_i$ in $\R^n$. Therefore we can apply Theorem~\ref{t:alg} to conclude that every 
homogeneous quasi-homomorphism vanishes on $\Ker\Cal$.

Applying the same argument to $G=\Diff_0^c(\R^n,\omega)$ with each $G_i$ equal to $\Diff_0^c(D_i,\omega)$, we conclude that the 
restriction to $G_i$ of any homogeneous quasi-ho\-mo\-mor\-phism $\varphi$ on $G$ is a homomorphism, and is therefore a constant multiple 
of the Calabi invariant. But any two balls in $\R^n$ are both contained in some larger ball, and therefore the restriction
of $\varphi$ is the same multiple of $\Cal$ on all balls.
\end{proof}
This argument clearly does not work in the finite volume case because the number of commuting conjugate embeddings 
one can construct is bounded in terms of the available volume. The abstract isomorphism type of the group of symplectic
or volume-preserving diffeomorphisms of a ball is always the same. In particular it does not depend on the size or volume 
of the ball. Therefore one can always find arbitrarily large numbers of commuting embeddings of such a group in the corresponding
group of a larger ball, or in itself, but these embeddings are conjugate only in the full diffeomorphism group, and not 
in the diffeomorphism group preserving $\omega$.

\section{Comparison with results on quasi-norms}\label{s:BIP}

In this section we compare our results on the non-existence of homogeneous quasi-ho\-mo\-mor\-phisms on various 
stable groups to results about existence and non-existence of quasi-norms. The following definitions are due to
Burago, Ivanov and Polterovich~\cite{BIP}.
\begin{definition}
A function $q\colon G\longrightarrow\R$ on a group is called a quasi-norm if it is almost subadditive, almost
invariant under conjugacy, and unbounded.

A group $G$ is called unbounded if it admits a quasi-norm, and is called bounded otherwise.
\end{definition}
Here almost subadditive and almost invariant under conjugation means that there is a constant $c$ such that
$$
q(xy)\leq q(x)+q(y)+c
$$
and 
$$
\vert q(xyx^{-1})-q(y)\vert\leq c
$$
for all $x$ and $y\in G$.

If $\varphi$ is a non-trivial homogeneous quasi-homomorphism, then its absolute value is a quasi-norm.
Therefore boundedness of a group is an even stronger property than the non-existence of homogeneous 
quasi-homomorphisms. A bounded group has finite Abelianization and vanishing stable commutator length.
If the Abelianization is trivial, then the group is uniformely perfect~\cite{BIP}.

Burago, Ivanov and Polterovich~\cite{BIP} prove boundedness of many diffeomorphism groups by an
argument that is rather stronger than ours. In particular, they prove boundedness for the groups we considered 
in Theorem~\ref{t:easy}. An easy example that illustrates the difference between the vanishing mechanism
in~\cite{BIP} and the one considered here is the following.
\begin{example}
Consider the group $G=\Diff_0^c(\R^n,\omega)$ from Theorem~\ref{t:Ham}. Then the function $q$ which 
assigns to every $g\in G$ the (symplectic) volume of its support is a quasi-norm which is non-trivial on
the kernel of the Calabi homomorphism. Therefore this kernel is an unbounded group which nevertheless 
has vanishing stable commutator length.

This also brings out another aspect of the difference between the finite and infinite volume cases mentioned earlier. 
If we take $G=\Diff_0^c(D,\omega)$ with $D$ of finite (symplectic) volume, then $q$ is bounded and therefore not
a quasi-norm.
\end{example}

This idea of a quasi-norm defined by support size actually applies to all the discrete groups we considered in Section~\ref{s:discrete}:
\begin{theorem}
The following groups are unbounded: 
\begin{enumerate}
\item the stable mapping class group $\Gamma_{\infty}$, 
\item the braid group $B_{\infty}$ on infinitely many strands, 
\item the stable automorphism group of free groups $Aut_{\infty}(F)$, and 
\item the stable special automorphism group of free groups $S\Aut_{\infty}(F)$.
\end{enumerate}
\end{theorem}
\begin{proof}
Perhaps the  case of the braid group is easiest to visualize. Define $q\colon B_{\infty}\longrightarrow\Z$ 
by sending a braid $x$ to the smallest number $k$ of strands needed to express it. These need not be the 
first $k$ strands, but can be any $k$ strands. The function $q$ is clearly invariant under conjugation and 
is unbounded. It is also subadditive, because a set of strands used to express $x$ and a set of strands
used to express $y$ can together be used to  express $xy$. Thus $q$ is a quasi-norm.

On the stable mapping class group define $q\colon\Gamma_{\infty}\longrightarrow\Z$ by mapping an element 
$x$ to the smallest genus $g$ of a compact surface with one boundary component on which a diffeomorphism representing 
the isotopy class $x$ can be supported. Again this does not have to be the embedding of $\Sigma^1_g$ from
the definition of the stabilization procedure, but can be any such subsurface of the infinite genus surface. 
By definition, this function is conjugacy-invariant. It is also not hard to see 
that it is unbounded, for example by considering the action of mapping classes on homology. To check (almost)
subadditivity, one has to consider several cases depending on how compact surfaces with one boundary
component supporting representatives for $x$ and $y$ sit in the infinite genus surface. In all cases the union
of these two subsurfaces can be enlarged slightly, without increasing the genus, to obtain a compact surface
with one boundary component supporting a representative for $xy$. The genus of this surface is at most
$q(x)+q(y)$. Therefore $q$ is in fact subadditive and a quasi-norm.

The argument for $Aut_{\infty}(F)$ and $S\Aut_{\infty}(F)$ is completely analogous.
\end{proof}

\bigskip
\noindent
{\bf Acknowledgements:}
I am very grateful to all the colleagues with whom I discussed the subject matter of this paper, including in particular
N.~A'Campo, D.~Burago, D.~Calegari, K.~Fujiwara, V.~Markovic, S.~Matsumoto, S.~Morita, L.~Polterovich, P.~Py and 
T.~Tsuboi. Special thanks are due to P.~Py for an inspiring conversation 
about the work of Burago, Ivanov and Polterovich, which got me started in the first place.

\bibliographystyle{amsplain}

\bigskip

\end{document}